\input amstex
 \documentstyle{amsppt}
\magnification 1200
\baselineskip 18pt
%%% \input definiti.tex
%
% file definiti.tex starts here 
%
\catcode`@=11
\def\binrel@#1{\setbox\z@\hbox{\thinmuskip0mu
\medmuskip\m@ne mu\thickmuskip\@ne mu$#1\m@th$}%
\setbox\@ne\hbox{\thinmuskip0mu\medmuskip\m@ne mu\thickmuskip
\@ne mu${}#1{}\m@th$}%
\setbox\tw@\hbox{\hskip\wd\@ne\hskip-\wd\z@}}
\def\overset#1\to#2{\binrel@{#2}\ifdim\wd\tw@<\z@
\mathbin{\mathop{\kern\z@#2}\limits^{#1}}\else\ifdim\wd\tw@>\z@
\mathrel{\mathop{\kern\z@#2}\limits^{#1}}\else
{\mathop{\kern\z@#2}\limits^{#1}}{}\fi\fi}
\def\underset#1\to#2{\binrel@{#2}\ifdim\wd\tw@<\z@
\mathbin{\mathop{\kern\z@#2}\limits_{#1}}\else\ifdim\wd\tw@>\z@
\mathrel{\mathop{\kern\z@#2}\limits_{#1}}\else
{\mathop{\kern\z@#2}\limits_{#1}}{}\fi\fi}
\def\circle#1{\leavevmode\setbox0=\hbox{h}\dimen@=\ht0
\advance\dimen@ by-1ex\rlap{\raise1.5\dimen@\hbox{\char'27}}#1}
\def\sqr#1#2{{\vcenter{\hrule height.#2pt
     \hbox{\vrule width.#2pt height#1pt \kern#1pt
       \vrule width.#2pt}
     \hrule height.#2pt}}}
\def\square{\mathchoice\sqr34\sqr34\sqr{2.1}3\sqr{1.5}3}
\def\force{\hbox{$\|\hskip-2pt\hbox{--}$\hskip2pt}}  
\catcode`@=\active

%
% file definiti.tex ends  here 
%
\define\pmf{\par\medpagebreak\flushpar}
\define\un{\underbar}

\define\pbf{\par\bigpagebreak\flushpar}
\define\a{\aleph}
\define\w{\omega}
\define\k{\kappa}
\topmatter
\title
Full Reflection of Stationary Sets Below $\a_\omega$
\endtitle
\author
Thomas Jech and Saharon Shelah
\endauthor
\thanks Supported in part by NSF and by a Fulbright grant at the Hebrew 
University.
Shelah's Publ. No. 387. Supported in part by the BSF.
\endthanks
\endtopmatter
\centerline{\un{Abstract}}

It is consistent that for every $n \ge 2$, every stationary subset of $\w_n$ 
consisting of ordinals of cofinality $\w_k$ where $k = 0$ or $k \le n -3$ 
reflects fully in the set of ordinals of cofinality $\w_{n-1}$.  
We also show that this result is best possible.
\par
\newpage
\pmf
1. \un{Introduction}.

A stationary subset $S$ of a regular uncountable cardinal $\k$ 
\un{reflects at} $\gamma < \k$ if 
$S \cap \gamma$ is a stationary subset of $\gamma$.  For stationary sets 
$S, A \subseteq \k$ let 
$$S < A \; \text{if} \; S \; \text{reflects at almost all} \;  \alpha \in A$$ 
where ``almost all" means modulo the closed 
unbounded filter on $\k$, i.e. with the exception of a nonstationary set 
of $\alpha$'s.  If $S < A$ we say that $S$ \un{reflects fully in} 
$A$.  The \un{trace} of $S, Tr (S)$, is the set of all $\gamma < \k$ at 
which $S$ reflects.  
The relation $<$ is well-founded [1], and $o (S)$, the \un{order of} $S$, 
is the rank of $S$ in this well-founded relation.

In this paper we investigate the question which stationary subsets of 
$\w_n$ reflect fully in which 
stationary sets; in other words the structure of the well founded relation 
$<$.  Clearly, $o (S) < o (A)$ 
is a necessary condition for $S < A$, and moreover, a set $S \subseteq 
\w_n$ has order $k$ just in case it has a stationary intersection with 
the set
$$ S^n_k = \{ \alpha < \w_n : cf \alpha = \w_k \}. $$ 
Thus the problem reduces to the investigation of full reflection 
of stationary subsets of $S^n_k$ in stationary subsets of $S^n_m$ for 
$k < m < n$.

The problem for $n =2$ is solved completely in Magidor's paper [2]:  It is 
consistent that every stationary $S \subseteq S^2_0$ reflects fully in 
$S^2_1$.  The problem for 
$n > 2$ is more complicated.  It is tempting to try the obvious 
generalization, namely $S < A$ 
whenever $o (S) < o (A)$, but this is provably false:
\pmf
\un{Proposition 1.1}.  There exist stationary sets $S \subset S^3_0$ and 
$A \subset S^3_1$ such that $S$ does not reflect at any $\gamma \in A$.
\pbf
\subheading{Proof}  Let $S_i, i < \w_2$, be any family of pairwise 
disjoint subsets of $S^3_0$, and let $\langle C_\gamma : \gamma \in 
S^3_1 \rangle$ be 
such that each $C_\gamma$ is a closed unbounded subset 
of $\gamma$ of order type $\w_1$.  Clearly, at most $\a_1$ of the sets $S_i$ 
can meet each $C_\gamma$, and so for each $\gamma$ there is 
$i (\gamma) < \w_2$ such that $C_\gamma \cap S_i = \emptyset$ for all 
$i \ge i (\gamma)$.  

There is $i < \w_2$ such that $i (\gamma) = i$ for a stationary set of 
$\gamma$'s.  Let $A \subset S^3_1$ be this stationary set and let 
$S = S_i$.  Then $S \cap C_\gamma = \emptyset$ for all $\gamma \in A$ 
and so $S \cap \gamma$ is nonstationary.  Hence $S$ does not reflect at 
any $\gamma \in A$.  \hfill $\square$

There is of course nothing special in the proof about $\a_3$ (or about 
$\a_1$) and so we have the following generalization:
\pmf
\un{Proposition 1.2}.  Let $k < m < n -1$.  There exist stationary sets 
$S \subseteq S^n_k$ and $A \subseteq S^n_m$ such that 
$S$ does not reflect at any $\gamma \in A$.  \hfill $\square$
\pmf
Consequently, if $n > 2$ then full reflection in $S^n_m$ is possible 
only if $m = n-1$.  This motivates our Main Theorem.
\pmf
1.3 \un{Main Theorem}.  Let $\k_2 < \k_3 < \cdots < \k_n < \cdots$ be a 
sequence of supercompact cardinals.  There is a 
generic extension $V [G]$ in which $\k_n = \a_n$ for all $n \ge 2$, and 
such that 

(a) every stationary subset of $S^2_0$ reflects fully in $S^2_1$, and 

(b) for every $n \ge 3$, every stationary subset of $S^n_k$ for all 
$k = 0, \cdots, n -3$, reflects fully in $S^n_{n-1}$.

We will show that the result of the Main Theorem is best possible.  But 
first we prove a corollary:
\pmf
\un{1.4 Corollary}.  In the model of the Main Theorem we have for all 
$n \ge 2$ and all $m, 0 < m < n$:

(a) Any $\a_m$ stationary subsets of $S^n_0$ reflect simultaneously at some 
$\gamma \in S^n_m$.

(b) For every $k \le m -2$,any $\a_m$ stationary subsets of $S^n_k$ reflect 
simultaneously at some $\gamma \in S^n_m$.
\pmf
\subheading{Proof}  Let us prove (a) as (b) is similar.  Let $m < n$ and let 
$S_\xi, \xi < \w_m$, be stationary subsets of $S^n_0$.  First, each 
$S_\xi$ reflects fully in 
$S^n_{n-1}$ and so there exist club sets $C_\xi, \xi < \w_m$, such that each 
$S_\xi$ reflects at all $\alpha \in C_\xi \cap S^n_{n-1}$.  As the club filter 
is $\w_n$ - complete, there exists an $\alpha \in S^n_{n-1}$ such that 
$S_\xi \cap \alpha$ is 
stationary, for all $\xi < \w_m$.  Next we apply full reflection of 
subsets of $S^{n-1}_0$ in 
$S^{n-1}_{n-2}$ (to the ordinal $\alpha$ of cofinality $\w_{n-1}$ rather 
than to 
$\w_{n-1}$ itself) and the $\w_{n-1}$ - completeness of the club filter on 
$\w_{n-1}$, to find 
$\beta \in S^n_{n-2}$ such that $S_\xi \cap \beta$ is stationary for all 
$\xi < \w_m$.  
This way we continue until we find a $\gamma \in S^n_m$ such that every 
$S_\xi \cap \gamma$ is stationary. \hfill $\square$

Note that the amount of simultaneous reflection in 1.4 is best possible:
\pmf
1.5 \un{Proposition}.  If $cf \gamma = \a_m$ and if $S_\xi, \xi < 
\w_{m+1}$, are disjoint 
stationary sets then some $S_\xi$ does not reflect at $\gamma$.
\pmf
\subheading{Proof}  $\gamma$ has a club subset of size $\a_m$, and it 
can only 
meet $\a_m$ of the sets $S_\xi \cap \gamma$.  \hfill $\square$

By Corollary 1.4, the model of the Main Theorem has the property that whenever 
$2 \le m < n$, every stationary subset of $S^n_k$ reflects quite strongly 
in $S^n_m$, provided $k \le m -2$.  
This cannot be improved to include the case of $k = m -1$, as the 
following proposition shows:
\pmf
1.6 \un{Proposition}.  Let $m \ge 2$.  Either

(a) for all $k < m -1$ there exists a stationary set $S \subseteq 
S^m_k$ that does not reflect fully in $S^m_{m-1}$, 

or

(b) for all $n > m$ there exists a stationary set $A \subseteq S^n_{m-1}$ 
that does not reflect at any $\delta \in S^n_m$.
\pmf
We shall give a proof of 1.6 in Section 3.  In our model we have, for 
every $m \ge 2$, full reflection 
of subsets of $S^m_0$ in $S^m_{m-1}$ (and of subsets of $S^m_k$ for $k 
\le m -3$) and therefore 1.6 (a) fails in the model.  Thus the model 
necessarily satisfies 1.6 (b), 
which shows that the consistency result is best possible. 
\pbf
2. \un{Proof of Main Theorem}

Let $\k_2 < \k_3  < \cdots < \k_n < \cdots$ be a sequence of cardinals 
with the 
property that for each $n \ge 2, \k_n$ is a $< \k_{n+1}$ - supercompact 
cardinal, i.e. for every 
$\gamma < \k_{n+1}$ there exists an elementary embedding $j : V \to M$ 
with critical point $\k_n$ 
such that $j (\k_n) > \gamma$ and $M^\gamma \subset M$. 
\footnote{We note in passing that the condition about the $\k_n$ is
equivalent to ``every $\k_n$ is $< \k_\w$ - supercompact" where
$\k_\w = sup_{m < \w} \k_m$.}
We construct 
the generic extension by iterated forcing, an iteration of length $\w$ 
with full support.  
The first stage of the iteration $P_1$ makes $\k_2 = \a_2$, and for each 
$n$, the $n^{th}$ stage $P_n$ (a forcing notion in $V(P_1 \ast \cdots 
\ast P_{n-1})$) makes $\k_{n+1} = \a_{n+1}$.  In the iteration, we 
repeatedly use three standard notions of forcing:  
$Col \; (\k, \alpha), C (\k)$ and $C U (\k, T)$.
\pmf
\un{Definition}.  Let $\k$ be a regular uncountable cardinal.

(a) $Col \;  (\k, \alpha)$ is the forcing that collapses 
$\alpha \ge \k$ with conditions of size $< \k$:  
\pmf
A condition is a function $p$ from a 
subset of $\k$ of size $< \k$ into $\alpha$; a condition $q$ is 
stronger than $p$ if $q \supseteq p$.

(b) $C (\k)$ is the forcing that adds a Cohen subset of $\k$:  A condition 
is an 0-1-function $p$ on a subset of $\k$ 
of size $< \k$; a condition $q$ is stronger than $p$ if $q \supseteq p$.

(c)  $C U (\k, T)$ is the forcing that shoots a club through a stationary set 
$T \subseteq \k$:

A condition is a closed bounded subset of $T$; a condition $q$ is stronger 
than $p$ if $q$ end-extends $p$.

The first stage $P_1$ of the iteration $P = \langle P_n: n = 1,2, \cdots 
\rangle$ is a forcing of size $\k_2$ that is $\w$ - closed
\footnote{A forcing notion is \un{$\lambda$ - closed} if every descending
sequence of length $\le \lambda$ has a lower bound.},
satisfies the 
$\k_2$ - chain condition and collapses each cardinal between $\a_1$ and 
$\k_2$ (it is essentially the Levy forcing with countable conditions.)  
For each $n \ge 2$, we construct (in $V (P | n)$) the $n^{th}$ stage 
$P_n$ such 
that
\pmf
(2.1) (a) $|P_n| = \k_{n+1}$

(b) $P_n$ is $\a_{n-2}$ closed

(c) $P_n$ satisfies the $\k_{n+1}$ - chain condition

(d) $P_n$ collapses each cardinal between $\a_n ( = \k_n)$ and $\k_{n+1}$

(e)  $P_n$ does not add any $\w_{n-1}$ - sequences of ordinals 
\pmf
and such that $P_n$ guarantees the reflection of stationary subsets of 
$\a_n$ stated in the theorem.

It follows, by induction, that each $\k_n$ becomes $\a_n$:  Assuming that 
$\k_n = \a_n$ in $V (P | n)$, the 
$n^{th}$ stage $P_n$ preserves $\a_n$ by (e), and the rest of the 
iteration $\langle P_{n+1}, P_{n+2}, \cdots \rangle$ also 
preserves $\a_n$ because it is $\a_{n-1}$ - closed by (b); $P_n$ makes 
$\k_{n+1}$ the successor of $\k_n$ by (c) and (d).

We first define the forcing $P_1$:  

$P_1$ is an iteration, with countable support, $\langle Q_\alpha: 
\alpha < \k_2 \rangle$ where for each $\alpha$, 
$$ Q_\alpha = Col \; (\a_1, \a_1 + \alpha) \times C (\a_1). $$  
It follows easily from well known facts that $P_1$ is an $\w$-closed 
forcing of size 
$\k_2$, satisfies the $\k_2$ - chain condition and makes $\k_2 = \a_2$.

Next we define the forcing $P_2$.  (It is a modification of Magidor's 
forcing from [2], but the added collapsing of cardinals requires a 
stronger assumption 
on $\k_2$ than weak compactness.  The iteration is padded up by the 
addition of Cohen forcing which 
will make the main argument of the proof work more smoothly).  
The definition of $P_2$ is inside the model $V (P_1)$, and so $\k_2 = \a_2$:

$P_2$ is an iteration, with $\a_1$ - support, $\langle Q_\alpha: \alpha 
< \k_3 \rangle$ where for each $\alpha$, 
$$ Q_\alpha = Col \; (\a_2, \a_2 + \alpha) \times C (\a_2) \times C U 
(T_\alpha) $$
where $T_\alpha$ is, in $V (P_1 \ast P_2 | \alpha)$, some stationary 
subset of $\w_2$.  We choose the $T_\alpha$'s so that 
each $T_\alpha$ contains all limit ordinals of cofinality $\w$.  It 
follows easily that for each $\alpha < \k_3, P_2 | \alpha \force 
Q_\alpha$ is $\w$-closed.

The crucial property of the forcing $P_2$ will be the following:
\pmf
\un{Lemma 2.2}.  $P_2$ does not add new $\w_1$ - sequences of ordinals.  

One 
consequence of Lemma 2.2 is that the conditions $(p,q,s) \in Q_\alpha$ 
can be taken to be sets in $V(P_1)$ (rather than 
in $V (P_1 \ast P_2 | \alpha)$).  Once we have Lemma 2.2, the 
properties (2.1) (a) - (e) follow easily.  

It remains to specify the choice of the $T_\alpha$'s.  By a standard 
argument using the $\k_3$ - chain condition, we can enumerate all 
potential subsets of $\w_2$ by a sequence $\langle S_\alpha : \alpha < 
\k_3 \rangle$ in such a way 
that each $S_\alpha$ is already in $V (P_1 \ast P_2 | \alpha)$.  At the 
stage $\alpha$ of the iteration, we let $T_\alpha = \w_2$ , unless 
$S_\alpha$ is, in $V (P_1 \ast P_2 | \alpha)$, a stationary set of 
ordinals of cofinality $\w$.  If that is the case, we let 
$$ T_\alpha = (Tr (S_\alpha) \cap S^2_1) \cup S^2_0 $$

Assuming that Lemma 2.2 holds, we now show that in $V (P_1 \ast P_2)$, 
every stationary $S \subseteq S^2_0$ reflects fully in $S^2_1$:  

The set $S$ appears as $S_\alpha$ at some stage $\alpha$, and because it is 
stationary in $V (P_1 \ast P_2)$, it is stationary in the smaller model 
$V (P_1 \ast P_2 | \alpha)$.  The forcing 
$Q_\alpha$ creates a closed unbounded set $C$ such that $C \cap S^2_1 
\subseteq Tr (S)$ (note that because $P_2$ does not 
add $\w_1$ - sequences, the meaning of $Tr (S)$ or of $S^2_1$ does not 
change). 

Thus in $V (P_1 \ast P_2)$ we have full reflection of subsets of $S^2_0$ 
in $S^2_1$.  The later stages of the 
iteration do not add new subsets of $\w_2$ and so this full reflection 
remains true in $V (P)$.

We postpone the proof of Lemma 2.2 until after the definition of the 
rest of the iteration.

We now define $P_n$ for $n \ge 3$.  We work in $V (P_1 \ast \cdots \ast 
P_{n-1})$.  By the induction hypothesis we have $\k_n = \a_n$.

$P_n$ is an iteration with $\a_{n-1}$ - support, $\langle Q_\alpha : 
\alpha < \k_{n+1} \rangle$, where for each $\alpha$, 
$$ Q_\alpha = Col (\a_n, \a_n + \alpha ) \times C (\a_n) \times C U 
(T_\alpha)$$
where $T_\alpha$ is a  $P_n | \alpha$ - name for a subset of $\w_n$.  
To specify the $T_\alpha$'s, let 
$\langle S_\alpha : \alpha < \k_{n+1} \rangle$ be an enumeration of all 
potential subsets of $\w_n$ such that each $S_\alpha$ 
is a $P_n | \alpha$ - name.  At the stage $\alpha$, let $T_\alpha = 
\w_n$ unless 
$S_\alpha$ a stationary set of ordinals and $S_\alpha \subseteq S^n_k$ 
for some $k = 0, \cdots, n -3$, in which case let 
$$ \align T_\alpha &= (Tr (S_\alpha) \cap S^n_{n-1}) \cup 
(S^n_0 \cup \cdots \cup S^n_{n-2}) \\
&= \{ \gamma < \w_n : cf \gamma \le \w_{n-2} \; \text{or} \; S_\alpha 
\cap \gamma \; \text{is stationary} \} \endalign $$
Due to the selection of the $T_\alpha$'s, $Q_\alpha$ is $\w_{n-2}$ - 
closed, and so is $P_n$.  
The crucial property of the forcing is the analog of Lemma 2.2:
\pmf
\un{Lemma 2.3}.  $P_n$ does not add new $\w_{n-1}$ - sequences of ordinals.

Given this lemma, properties (2.1) (a) - (e) follow easily.  The same 
argument as 
given above for $P_2$ shows that in $V (P_1 \ast \cdots \ast P_n)$, and 
therefore in $V(P)$ as well, every 
stationary subset of $S^n_k, k = 0, \cdots, n -3$, reflects fully in 
$S^n_{n-1}$.

It remains to prove Lemmas 2.2  and 2.3.  We prove Lemma 2.3, as 2.2 is 
an easy modification.
\pbf
\subheading{Proof of Lemma 2.3}

Let $n \ge 3$, and let us give the argument for a specific $n$, say $n =4$.  
We want to show that 
$P_4$ does not add $\w_3$ -sequences of ordinals.  

We will work in $V (P_1 \ast P_2)$ (and so consider the forcing $P_3 
\ast P_4)$.  As $P_1 \ast P_2$ has size $\k_3, \k_4$ 
is a $< \k_5$ - supercompact cardinal in $V (P_1 \ast P_2)$, and 
$\k_3 = \a_3$.  The 
forcing $P_3$ is an iteration of length $\k_4$ that makes $\k_4 = \a_4$ 
and is $\a_1$ - closed; then $P_4$ is an iteration of length $\k_5$.  By 
induction on $\alpha < \k_5$ we show  
$$ P_4 | \alpha\; \text{ does not add} \;  \w_3 \; \text{- sequences 
of ordinals.} \tag 2.4 $$  As $P_4$ has the $\a_5$ - chain condition, 
(2.4) is certainly enough for Lemma 2.3.  Let $\alpha < \k_5$.

Let $j$ be an elementary embedding $j: V \to M$ (as we work in $V (P_1 
\ast P_2), V$ means 
$V (P_1 \ast P_2)$) such that $j (\k_4) > \beta$ and $M^\beta \subset M$,  
for some inaccessible cardinal $\beta > \alpha$.  
Consider the forcing $j (P_3)$ in $M$.  It is an iteration of which 
$P_3$ is an initial segment.  By a standard argument, 
the elementary embedding $j : V \to M$ can be extended to an elementary 
embedding $j: V (P_3) \to M (j (P_3))$.  We claim that every 
$\beta$-sequence of ordinals in $V (P_3)$ belongs to $M (j (P_3)$): 
the name for such a set has size $\le \beta$ 
and so it belongs to $M$, and since $P_3 \in M$ 
and $M (P_3) \subseteq M (j (P_3))$, the claim follows.  In particular, 
$P_4 | \alpha \in M (j (P_3))$.

Let $p,\dot F \in V (P_3)$ be such that $p \in P_4 | \alpha$ and 
$\dot F$ is a $(P_4 | \alpha)$ - name 
for an $\w_3$ - sequence of ordinals.  We shall find a stronger condition 
that decides all the values of $\dot F$.  By the 
elementarity of $j$, it suffices to prove that 

(2.5) $\exists \bar p \le j (p)$ in $j (P_4 | \alpha)$ that decides 
$j (\dot F)$.  
\pmf
The rest of the proof is devoted to the proof of (2.5).  

Let $G$ be an $M$ - generic filter on $j (P_3)$.
\pmf
\un{Lemma 2.6}.  In $M [G]$ there is a generic filter $H$ on $P_4 | \alpha$ 
over $M [G \cap P_3]$ such that $M [G]$ is a generic extension of 
$M [G \cap P_3] [H]$ by an $\a_1$ - closed forcing, and such that $p \in H$.
\pbf
\subheading{Proof}  There is an $\eta < j (\k_4)$ such that $P_4 | \alpha$ 
has size $\a_3$ in $M_\eta = M [G \cap (j (P_3) | \eta ) ]$.  Since 
$P_4 | \alpha$ is $\a_2$ - closed, it is isomorphic in $M_\eta$ to the 
Cohen forcing $C (\a_3)$.  But $Q_\eta = (j (P_3)) (\eta) = Col (\a_3, 
\a_3 + \eta) \times C (\a_3) \times C U (T_\eta)$, so 
$G | Q_\eta = G_{Col} \times G_C \times G_{C U}$, and using $G_C$ and
the isomorphism between $P_4 | \alpha$ and $C (\a_3)$ we obtain $H$.  
Since the quotient forcing $j (P_3) / 
(P_3 \times C (\a_3))$ is an iteration of $\a_1$ - closed forcings, it is
$\a_1$ - closed. \hfill $\square$
\pmf
\un{Lemma 2.7}.  In $M [G]$ there is a condition $\bar p \in j (P_4 | 
\alpha)$ that extends $p$, and extends every member of $j^{''} H$.  

Lemma 2.7 will complete the proof of (2.5):  since every value of $\dot F$ 
is decided by some condition in $H$, every value of 
$j (\dot F)$ is decided by some condition in $j^{''} H$, and therefore by 
$\bar p$.
\pbf
\subheading{Proof of Lemma 2.7}  Working in $M [G]$, we construct $\bar p 
\in j (P_4 | \alpha)$, a 
sequence $\langle p_\xi : \xi < j (\alpha) \rangle$ of length $j (\alpha)$, 
by induction.  When $\xi$ is not in the range of $j$, we let 
$p_\xi$ be the trivial condition; that guarantees that the support of 
$\bar p$ has size $| \alpha|$ which is $\a_3$ (because $\alpha < j 
(\k_4) = \a_4$ 
in $M [G]$).  So let $\xi < \alpha$ be such that 
$\bar p | j (\xi)$ has been defined, and construct $p_{j (\xi)}$.  

The condition $p_{j (\xi)}$ has three parts $u, v, s$ where $u \in 
Col (j (\k_4), j (\k_4) + j (\xi)), v \in C ( (\k_4))$ and $s \in 
C U (T_{j (\xi)})$.  It is easy to construct the $u$ - 
part and the $v$ - part, as follows: The filter $H | P_4 (\xi)$ has three
parts; a collapsing function $f$ of $\k_4$ onto $\k_4 + \xi$, a 0-1-function
$g$ on $\k_4$, and a club subset $C$ of $T_\xi$. We let $u = j " f$ and
$v = j " g$, and these
are functions of size $\a_3$ 
and therefore members of $Col$ and $C$ 
respectively.  For the $s$ - part, let $s = j " C \cup \{ \k_4 \}$.
In order that this set be a condition in 
$C U (T_{j (\xi)} )$, we have to verify 
that $\k_4 \in T_{j (\xi)}$.

This is a nontrivial requirement if $S_{j (\xi)}$ is in 
$M (j (P_3) \ast (j (P_4) | j (\xi)))$ a stationary subset of 
$j (\k_4)$ and is a subset of 
either $S^4_0$ or of $S^4_1$ (of $S^n_k$ for $n =4$ and $k \le n -3$).  
Then $\k_4$ has to be reflecting point of $S_{j (\xi)}$, 
i.e. we have to show that $S_{j (\xi)} \cap \k_4$ 
is stationary, in $M (j (P_3) \ast (j (P_4) | j (\xi))$.

By the assumption and by elementarity of $j, S_\xi$ is a stationary 
subset of $\k_4$ in $V (P_3 \ast P_4 | \xi)$, and 
$S_\xi \subseteq S^4_0$ or $S_\xi \subseteq S^4_1$, i.e. consists of 
ordinals of cofinality $\le \w_1$.  Since $S_{j (\xi)} \cap \k_4 = 
j (S_\xi) \cap \k_4 = S_\xi$, 
it suffices to show that $S_\xi$ is stationary not only in 
$V (P_3 \ast P_4 | \xi)$ but also in $M (j (P_3) \ast (j (P_4) | j (\xi))$.

Firstly $M (P_3 \ast P_4 | \xi) \subseteq V (P_3 \ast P_4 | \xi)$, and so 
$S_\xi$ is stationary in $M (P_3 \ast P_4 | \xi)$.  Secondly, $j (P_4)$ is 
$\a_1$ - closed, and by Lemma 2.6, $M (j (P_3))$ is an $\a_1$ - closed 
forcing extension of $M (P_3 \ast P_4 | \xi)$, and so the proof is 
completed by application of the following lemma (taking $\k = \a_0$ or 
$\a_1, \lambda = \a_4$).
\pmf
\un{Lemma 2.8}  Let $\k < \lambda$ be regular cardinals and assume that 
for all $\alpha < \lambda$ and all $\beta < \k, \alpha^\beta < \lambda$.  
Let $Q$ be a $\k$ - closed forcing and $S$ a stationary subset of 
$\lambda$ of ordinals of 
cofinality $\k$.  Then $Q \force S$ is stationary.

This lemma is due to Baumgartner and we include the proof for lack of 
reference.
\pbf
\subheading{Proof of Lemma 2.8}  Let $q$ be a condition and let $\dot C$ 
be a $Q$ - name for a closed unbounded subset of $\lambda$.  
We shall find $\bar q \le q$ and $\gamma \in S$ such that $\bar q \force 
\gamma \in \dot C$.  Let $M$ be a 
transitive set such that $M$ is a model of enough set theory, is closed 
under $< \k$ - sequences and such 
that $M \supseteq \lambda, q \in M, Q \in M, \dot C \in M$.  Let 
$\langle N_\gamma : \gamma < \lambda \rangle$ be an 
elementary chain of submodels of $M$ such that each $N_\gamma$ has size 
$< \lambda$, contains $q, Q$ and $\dot C, N_\gamma \cap \lambda$ is an 
ordinal, and $N_{\gamma +1}$ contains all $< \k$ - sequences in $N_\gamma$.  
Since $S$ is stationary, there exists a $\gamma \in S$ such 
that $N_\gamma \cap \lambda = \gamma$.  As $cf \gamma = \k, N = N_\gamma$ 
is closed under $< \k$ - sequences.       

\par
Let $\{ \gamma_\xi : \xi < \k \}$ be an increasing sequence with limit 
$\gamma$.  We construct 
a descending sequence $\{ q_\xi : \xi < \k \}$ of conditions such that 
$q_0 =  q$, such that for all $\xi < \k, q_\xi \in N$ and for some 
$\beta_\xi \in N$ greater than 
$\gamma_\xi , q_{\xi +1} \force \beta_\xi \in \dot C$.  At successor 
stages, $q_{\xi +1}$ exists because in $N, q_\xi$ forces 
that $\dot C$ is unbounded.  At limit stages $\eta < \k$, the $\eta$ - 
sequence $\langle q_\xi : \xi < \eta \rangle$ is in $N$ and has a lower 
bound in $N$ because $N \models Q$ is $\k$ - closed.

Since $Q$ is $\k$ - closed, the sequence $\langle q_\xi : \xi<\k \rangle$ 
has a 
lower bound $\bar q$, and because of the $\beta$'s, $\bar q$ forces that 
$\dot C$ is unbounded in $\gamma$.  
Therefore $\bar q \force \gamma \in \dot C$.  \hfill $\square$
\pbf
3. \un{Negative results}.

We shall now present several negative 
results on the structure of the relation $S < T$ below $\a_\w$.  With 
the exception of the proof of Proposition 1.6, we state the results for 
the particular case of 
reflection of subsets of $S^3_0$ in $S^3_1$, but the results generalize 
easily to other cardinalities and other cofinalities.

The first result uses a simple calculation (as in Proposition 1.1):
\pmf
\un{Proposition 3.1}.  For any $\a_3$ stationary sets $A_\alpha \subseteq 
S^3_1, \alpha < \w_3$, 
there exists a stationary set $S \subseteq S^3_0$ such that $S \not< 
A_\alpha$ for all $\alpha$.
\pbf
\subheading{Proof}  Let $A_\alpha, \alpha < \w_3$, be stationary subsets 
of $S^3_1$.  By [3], there 
exist $\a_4$ almost disjoint stationary subsets of $S^3_0$; let $S_i, i 
< \w_4$, be such sets.  
Assuming that each $S_i$ reflects fully in some $A_{\alpha (i)}$, we can 
find $\a_4$ of them that reflect fully in the same $A_\alpha$.  
Take any $\a_2$ of them and reduce each by a nonstationary set to get 
$\a_2$ pairwise disjoint stationary subsets $\{ T_\xi : \xi < \w_2 \}$ 
of $S^3_0$, such that 
each of them reflects fully in $A_\alpha$.  Hence there are clubs $C_\xi, 
\xi < \w_2$, such that $Tr (T_\xi) \supseteq A_\alpha \cap C_\xi$ for 
every $\xi$.  Let 
$\gamma \in \underset \xi < \w_2 \to \cap C_\xi \cap A_\alpha$.  Then 
every $T_\xi$ reflects at $\gamma$, and so $\gamma$ has $\a_2$ pairwise 
disjoint stationary subsets $\{ T_\xi \cap \gamma: \xi < \w_2 \}$.  This 
is a contradiction 
because $\gamma$ has a closed unbounded subset of size $cf \gamma = 
\a_1$.  \hfill $\square$

The next result uses the fact that under GCH there exists a 
$\diamondsuit$ - sequence for $S^3_1$.
\pmf
\un{Proposition 3.2}.  (GCH)  There exists a stationary set $A \subseteq 
S^3_1$ that is not the trace of any $S \in S^3_0$; precisely: for every 
$S \subseteq S^3_0$ the set $A \Delta (Tr (S) \cap S^3_1)$ is stationary.
\pbf
\subheading{Proof}  Let $\langle S_\gamma : \gamma \in S^3_1 \rangle$ 
be a $\diamondsuit$ - sequence for 
$S^3_1$; it has the property that for every set $S \subseteq \w_3$, the 
set $D (S) = \{ \gamma \in S^3_1: S \cap \gamma = S_\gamma \}$ is stationary.  Let 
$$ A = \{ \gamma \in S^3_1 : S_\gamma \; \text{is nonstationary} \}. $$
The set $A$ is stationary because $A \supseteq D (\emptyset)$.  If $S$ is any 
stationary subset of $S^3_0$, then for every $\gamma$ in the stationary 
set $D (S), \gamma \in A$ iff $\gamma \notin Tr (S)$, and so 
$D (S) \subseteq A \Delta Tr (S)$. \hfill $\square$

The remaining negative results use the following theorem of Shelah which 
proves the existence of sets with the ``square 
property".
\pmf
\un{Theorem} ([4], Lemma 4.2).  Let $1 \le k \le n -2$.  The set $S^n_k$ is the 
union of $\a_{n-1}$ stationary sets $A$, each having the following 
property.  There exists a collection $\{ C_\gamma : \gamma \in A \}$ 
(a ``square sequence for $A$") 
such that for each $\gamma \in A, C_\gamma$ is a club subset of 
$\gamma$ of order type $\w_k$, consisting of limit ordinals of cofinality 
$< \w_k$, and such that for all $\gamma_1, \gamma_2 \in A$ and all 
$\alpha$, if $\alpha \in C_{\gamma_1} \cap C_{\gamma_2}$ then 
$C_{\gamma_1} \cap \alpha = C_{\gamma_2} \cap \alpha$.

Square sequences can be used to construct a number of counterexamples.  
For instance, if $S_n, n < \w$, are $\a_0$ stationary subsets of $S^3_0$, 
then $Tr (\overset \infty \to{\underset n =0 \to \bigcup} S_n ) = 
\overset \infty 
\to {\underset n =0 \to \bigcup} S_n$.  Using a square sequence we get:
\pmf
\un{Proposition 3.3}.  There is a stationary set $A \subseteq S^3_1$ and 
stationary subsets $S_i, i < \w_1$, of $S^3_0$ such that $Tr (S_i) \cap A 
= \emptyset$ for each $i$ but $Tr (\underset i < \w_1 \to \bigcup S_i) 
\supseteq A$.
\pbf
\subheading{Proof}  Let $A$ be a stationary subset of $S^3_1$ with a square 
sequence $\{ C_\gamma : \gamma \in A \}$, and let $S = \underset \gamma 
\in A \to \bigcup C_\gamma$.  Clearly, $S \subseteq S^3_0$ is stationary, and 
$Tr (S) \supseteq A$.  For each $\xi < \w_1$, let 
$$ S_\xi = \{ \alpha \in S: \; \text{order type} \; (C_\gamma \cap 
\alpha) = \xi \} $$
(this is independent of the choice of $\gamma \in A$).  For every 
$\gamma \in S$ and every $\xi < \w_1$, the set $S_\xi \cap C_\gamma$ has 
exactly one element, and 
so $S_\xi$ does not reflect at $\gamma$.    It is easy to see that $\a_1$ 
of the sets $S_\xi$ are stationary.  [The definition of $S_\xi$ is a 
well known trick] \hfill $\square$

The argument used in the above proof establishes the following:
\pmf
\un{Proposition 3.4}.  If a stationary set $A \subseteq S^n_m$ has a square 
sequence and if $k < m$ then 
there exists a stationary $S \subseteq S^n_k$ that does not reflect at 
any $\gamma \in A$. \hfill $\square$
\pbf
\subheading{Proof of Proposition 1.6}  Let $2 \le m <n$ and let us 
assume that (b) fails, i.e. that every stationary set $A \subseteq 
S^n_{m-1}$ reflects at some $\delta$ of cofinality $\a_m$.  We shall 
prove that (a) holds.  For each $k < m -1$ we want a stationary set 
$S \subseteq 
S^m_k$ that does not reflect fully in $S^m_{m-1}$.  Let $k < m -1$.

Let $A$ be a stationary subset of $S^n_{m-1}$ that have a square sequence 
$\{ C_\gamma : \gamma \in A \}$.  
The set $A$ reflects at some $\delta$ of cofinality $\w_m$.  Let $C$ be 
a club subset of $\delta$ of order 
type $\w_m$.  Using the isomorphism between $C$ and $\w_m$, the sequence 
$\{C_\gamma \cap 
C: \gamma \in A \}$ becomes a square sequence for a stationary subset 
$B$ of $S^m_{m-1}$.  It follows that 
there is a stationary subset of $S^m_k$ that does not reflect at any 
$\gamma \in B$. \hfill $\square$ 

The last counterexample also uses a square sequence.   

\un{Proposition 3.5}.  (GCH)  There is a stationary set $A \subseteq 
S^3_1$ and $\a_4$ stationary sets $S_i \subseteq S^3_0$ such that the sets 
$\{ Tr (S_i) \cap A: i < \w_4 \}$ are stationary and pairwise almost disjoint.
\pbf
\subheading{Proof}  Let $A$ be a stationary subset of $S^3_1$ with a 
square sequence $\langle C_\gamma : \gamma \in 
A \rangle$, and let $S = \underset \gamma \in A \to \bigcup C_\gamma$.  
Let $\{ f_i : i < \w_4 \}$ be regressive functions on $S^3_0 \cup S^3_1$ 
with the property that for any two $f_i, f_j$, the set 
$\{ \alpha : f_i (\alpha) = f_j (\alpha) \}$ is nonstationary (such a 
family exists by [3]).  For each $i$ and each 
$\gamma \in A$, the function $f_i$ is regressive on $C_\gamma$ and so there 
is some $\eta = \eta (i, \gamma) < \gamma$ such that 
$\{ \alpha \in C_\gamma: f_i (\alpha) < \eta \}$ is stationary.  Let 
$T_{i, \gamma} \subseteq \w_1$ be the stationary set $\{ o.t. 
(C_\gamma \cap \alpha) : f_i (\alpha) < \eta \}$ and 
let $H_{i, \gamma}$ be the function on $T_{i, \gamma}$ (with values $< \eta$) 
defined by $H (\xi) = f_i (\xi^{th}$ element of $C_\gamma$).  For 
each $i$, the function on $A$ that to each $\gamma$ assigns $(T_{i 
\gamma}, H_{i \gamma})$ is regressive, and so constant 
$= (T_i, H_i)$ on a stationary set.  By a counting argument, $(T_i, H_i)$ is 
the same for $\a_4 \; i$'s; so w.l.o.g. 
we assume that they are the same $(T, H)$ for all $i$.

Now we let, for each $i$, 
$$ \align A_i &= \{ \gamma \in A: (\forall \alpha \in C_\gamma) \; \text{if} 
\; \xi = o.t. (C_\gamma \cap \alpha) \in T \; \text{then} \; 
f_i (\alpha) = H (\xi) \} \\
S_i &= \{ \alpha \in S: o.t. (C_\gamma \cap \alpha) \in T \; \text{and} 
(\forall \beta \le \alpha, \beta \in C_\gamma) \; \text{if} \; \xi = o.t. 
(C_\gamma \cap \beta) \in T \; \text{then} \; f_i (\beta) = 
H (\xi) \} \endalign $$
By the definition of $T$ and $H$, each $A_i$ is a stationary set, and 
each $S_i$ reflects at every point of 
$A_i$.  We claim that if $\gamma \in A$ and $S_i \cap \gamma$ is 
stationary then $\gamma \in A_i$. So let $\gamma \in A$ be such that 
$S_i \cap \gamma$ is stationary.  Let $\xi \in T$ and let 
$\alpha$ be the $\xi^{th}$ element 
of $C_\gamma$; we need to show that $f_i (\alpha) = H (\xi)$.  As $S_i 
\cap \gamma$ is stationary, 
there exists a $\beta \in S_i \cap C_\gamma$ greater than $\alpha$.  By the 
definition of $S_i, f_i (\alpha) = H (\xi)$.  
Thus $\gamma \in A_i$, and $A_i = A \cap Tr (S_i)$.  

Finally, we show that the sets $A_i$ are pairwise almost 
disjoint.  Let $C$ be a 
club disjoint from the set $\{ \alpha: f_i (\alpha) = f_j (\alpha) \}$. 
We claim that the set $C'$ of all limit points of $C$ is disjoint from 
$A_i \cap A_j$.  If $\gamma \in C'$ then 
$C \cap \gamma$ is a club in $\gamma$, and so is $C \cap C_\gamma$.  
Since $T$ is stationary in $\w_1$, there is a 
$\xi \in T$ such that the $\xi^{th}$ element $\alpha$ of $C_\gamma$ 
is in $C$, and therefore $f_i (\alpha) \neq f_j (\alpha)$; it follows 
that $\gamma$ cannot be 
both in $A_i$ and in $A_j$.  \hfill $\square$
\par
\newpage
\item{[1]}  T. Jech, Stationary subsets of inaccessible cardinals, in 
``Axiomatic Set Theory" (J. Baumgartner, ed.), Contemporary Math. 31, 
Amer. Math. Soc. 1984, 115-142.
\item{[2]}  M. Magidor, Reflecting stationary sets, J. 
Symbolic Logic \un{47} (1982), 755-771.
\item{[3]}  S. Shelah [Sh 247], More on stationary coding, in ``Around
Classification Theory of Models", Springer-Verlag Lecture 
Notes 1182 (1986), pp. 224-246.
\item{[4]}  S. Shelah [Sh 351], Reflecting stationary sets and 
successors of singular cardinals.     
\pbf
The Pennsylvania State University   \hfill The Hebrew University 
\bye